\documentclass[12pt]{amsart}
\usepackage{amsmath}
\usepackage{amsfonts,amssymb,amsthm,amstext,amscd,boxedminipage}  %
\usepackage[mathscr]{eucal}

\usepackage{graphicx,psfrag} %
\usepackage{color}
 \numberwithin{equation}{section}

\newcommand{\cross}{{\times}} %
\newcommand{\shat}{{\hat{s}}}
\renewcommand{\setminus}{{\smallsetminus}}  %
\DeclareMathSymbol{\minus} {\mathord}{operators}{"2D} %

\theoremstyle{plain}
\newtheorem{thm}{Theorem}[section]
\newtheorem{lem}[thm]{Lemma}
\newtheorem{cor}[thm]{Corollary}

\theoremstyle{definition}
\newtheorem{df}[thm]{Definition}
\newtheorem{remark}[thm]{Remark}
\newtheorem{example}[thm]{Example}

\input epsf

\def \H {\mathbb{H}}
\def \Z {\mathbb{Z}}
\def \D {\mathbb{D}}

\def \BRT {{Bollob\'as--Riordan--Tutte }}

\begin{document}

\title{The Jones polynomial and graphs on surfaces}
\date{\today}

\author[O. Dasbach]{Oliver T. Dasbach}
\address{Department of Mathematics, Louisiana State University,
Baton Rouge, LA 70803}
\email{kasten@math.lsu.edu}
\thanks {The first author was supported in part by NSF grants DMS-0306774 and DMS-0456275 (FRG)}

\author[D. Futer]{David Futer}
\address{Department of Mathematics, Michigan State University,
East Lansing, MI 48824}
\email{dfuter@math.msu.edu}
\thanks {The second author was supported in part by NSF grant DMS-0353717 (RTG)}

\author[E. Kalfagianni]{Efstratia Kalfagianni}
\address{Department of Mathematics, Michigan State University,
East Lansing, MI 48824}
\email{kalfagia@math.msu.edu}
\thanks{The third author was supported in part by NSF grants DMS-0306995 and DMS-0456155 (FRG)}

\author[X.-S. Lin]{Xiao-Song Lin}
\footnote{We regretfully inform you that Xiao-Song Lin  passed away on the 14th of January, 2007.}

\author[N. Stoltzfus]{Neal W. Stoltzfus}
\address{Department of Mathematics, Louisiana State University,
Baton Rouge, LA 70803}
\email{stoltz@math.lsu.edu}
\thanks {The fifth author was supported in part by NSF grant DMS-0456275 (FRG)}

\begin{abstract}
The Jones polynomial of an alternating link is a certain specialization of
the Tutte polynomial of the (planar) checkerboard graph associated to an
alternating projection of the link.  The \BRT polynomial generalizes the Tutte
polynomial of  graphs to graphs that are
embedded in closed oriented surfaces of higher genus. 

In this paper we show that the Jones polynomial of any link can be obtained from the \BRT 
polynomial of a certain oriented ribbon graph associated to a link projection. We give some applications of this approach. 
\end{abstract}

\maketitle

\section{Introduction}
 
Informally, oriented ribbon graphs (a.k.a. combinatorial maps, rotation systems) are graphs with a cyclic
orientation of the edges meeting at a vertex (see Chapter 10 of Tutte's book on graph theory \cite{Tutte:Book}). Their genus is the minimal genus of an oriented surface constructed by attaching polygonal faces in a manner prescribed by the cyclic orientation. Recently, Bollob\'as and Riordan extended the Tutte
polynomial to ribbon graphs \cite{BollobasRiordan:CyclicGraphs} as a three-variable polynomial,
where the exponents of the third variable are related to the genus of the ribbon graph. 

The ribbon graphs of interest in this introductory paper arise naturally from link projections.
We will show that the Jones polynomial, via the Kauffman bracket, is a specialization of the \BRT
polynomial of these oriented ribbon graphs. Furthermore, this approach leads to a natural notion of the minimal genus of the ribbon graph over all link projections 
which we call \emph{Turaev genus} of a link. Links of Turaev genus zero are exactly the alternating links. 

Our approach is different from the one taken by Thistlethwaite
\cite{Thistlethwaite:SpanningTreeExpansion}, where he gives a spanning tree
expansion of the Jones polynomial via signed graphs. Using these ideas, Kauffman
defined a three-variable Tutte polynomial for signed graphs
\cite{Kauffman:TutteSignedGraphs}.  For alternating links, our approach coincides with
Thistlethwaite's approach: the Jones polynomial of an alternating link
is a specialization of the Tutte polynomial of Tait's checkerboard graph of an
alternating link projection. The connection between the Tutte polynomial and the 
Jones polynomial for alternating knots was fruitfully used in \cite{DL:VolumeIsh, DasbachLin:HeadAndTail}.
The books \cite{Bollobas:Book,
Welsh:BookComplexity} give a good introduction to the interplay between knots
and graphs.   

There is a version of the Jones polynomial for links in 3-manifolds $M$
that are $I$-bundles over orientable surfaces: that is,  $M=S\times I$. In this
setting, one defines a version of the Kauffman bracket for link projections
on $S$ and the Jones polynomial is a normalization of this bracket that
remains invariant under Reidemeister moves on $S$. Chmutov and
Pak \cite{CS04} were the first to relate the \BRT graph polynomial to link invariants:
they showed that the Kauffman bracket of an alternating
link projection on $S$ is an evaluation of the \BRT polynomial of the
checkerboard graph of the projection. Thus, Chmutov and Pak generalize
Thistlethwaite's result for alternating links in $I$-bundles. In \cite{ChmutovPak:Virtual}
they extend their result to virtual links.

In this paper we relate the \BRT  polynomial and the Kauffman bracket in
a context different than that of Chmutov and Pak: we work with links in $S^3$ and consider 
projections up to the usual Reidemeister moves. We show that the Kauffman bracket of any
connected link projection is obtained as an evaluation of the \BRT polynomial
of a certain oriented ribbon graph (the $A$-graph) associated to the projection. We
conclude that the Jones polynomial of any link can be obtained from the \BRT
polynomial of the ribbon graph. From this point of view, our result is a
generalization of Thistlethwaite's alternating link result to all links.
 
The paper is organized as follows: Section \ref{sec:dessin} recalls the definition of an oriented ribbon graph. 
Section \ref{A-dessin} describes the construction of a ribbon graph from a diagram of a link. 
Section \ref{sec:duality} uses the surface of the ribbon graph to prove a duality result.  The relationship to Fomenko's concept of \textit{atom} \cite{F1991}, as introduced into knot theory by Manturov \cite{M98}, is also discussed.  In Section \ref{sec:BRT}, we define the \BRT polynomial of this ribbon graph, and show how the Kauffman bracket of the diagram can be obtained as a specialization of this polynomial. Passing from the Kauffman bracket to the Jones polynomial is then a matter of multiplying by a well-known diagrammatic factor. 

As an application, we get a spanning tree and a spanning subgraph expansion for the Kauffman bracket in Section \ref{sec:SpanningTrees}. Finally, Section \ref{sec:Adequate} gives some implications for
adequate links.

{\bf Acknowledgment: } We would like to thank Sergei Chmutov and Igor Pak for helpful discussions
on the \BRT polynomial. We also thank James Oxley and Vassily Manturov for helpful comments.

\section{Ribbon Graphs} \label{sec:dessin}

Heuristically, an \textit{oriented ribbon graph} can be viewed
as a multi-graph (i.e. loops and multiple edges are allowed) equipped with a
cyclic order on the edges at every vertex. Isomorphisms between oriented ribbon graphs are
graph isomorphisms that preserve the given cyclic order of the edges. 

 The definition of a ribbon graph, however, highlights the permutation (and
secondarily, the topological) structure. \begin{df} A \textit{connected oriented ribbon graph}
is a triple, $\D=(\sigma_0,\sigma_1,\sigma_2)$ of permutations of a finite set
$\mathcal{B} = [2n]:=\{1,2,\ldots,2n-1,2n\}$.  The triple must satisfy:

\begin{itemize}
\item $\sigma_1$ is a fixed point free involution, i.e. $\sigma_1(\sigma_1(b)) = b, \sigma_1(b)\ne b$ for all $b$.
\item $\sigma_0(\sigma_1(\sigma_2(b)))=b$ and
\item The group generated by $\langle \sigma_0, \sigma_1 \rangle$ acts transitively on B.
\end{itemize}
A \textit{oriented ribbon graph} is a disjoint union of connected oriented ribbon graphs. We will often use the term \emph{ribbon graph} for short, keeping the orientation implicit.
\end{df}

The elements of the set $\mathcal{B}$ will be called \textit{half edges}
(\textit{les brins} en fran{\c{c}}ais). Note that it follows from the second
condition that any two permutations determine the third.

\subsection{Associated Surface and Genus of a Ribbon Graph}

Given an oriented ribbon graph $\D=(\sigma_0,\sigma_1,\sigma_2)$,
the orbits of $\sigma_0$ form the vertex set, of cardinality $v(\D)$, the
orbits of $\sigma_1$, the edge set, of cardinality $e(\D)$ and $\sigma_2$, the
face set, of cardinality $f(\D)$, respectively. The underlying graph of a
ribbon graph has an edge connecting the vertices in whose orbit its two \textit{half
edges} lie. In addition, this graph is embedded in an oriented surface with
(cellular) faces corresponding to the orbits of $\sigma_2$ and oriented so that
$\sigma_0$ cyclically rotates the \textit{half edges} meeting at a vertex in
the rotation direction determined by the orientation.  By convention, we also allow the ribbon graph with an empty set of half-edges. 
The underlying graph has one vertex and one face.

\begin{df}\label{def:genus} The genus $g(\D)$ of a ribbon graph $\D$ with $k$ components is determined
by its Euler characteristic: $v(\D)-e(\D) +f(\D) = 2k-2g(\D)$. 
\end{df}

For the definitions and properties of oriented ribbon graphs, we will follow the recent
monograph \textit{Graphs in surfaces and their applications}
\cite{Lando:GraphSurface}. In particular, note that we will always assume the
embedding surface is oriented (certain extensions to the non-orientable case
are discussed in \cite{BollobasRiordan:NonOrientableSurfaces}). Historically,
this concept has been rediscovered several times and explored in several distinct
settings: combinatorics, topology, physics of fields, combinatorial group
theory and algebraic number theory. 

Oriented ribbon graphs have also 
been called:  fat graphs, cyclic graphs, combinatorial maps, rotation
systems and dessins d'enfant (Grothendieck's
terminology in the context of the Galois theory of certain arithmetic algebraic curves). In addition to 
the book of S.K. Lando and A. Zvonkin \cite{Lando:GraphSurface} and Tutte's book \cite{Tutte:Book}, surveys of note are those by Robert Cori and
Antonio Machi \cite{CoriMachi:Maps} and the book edited by Leila Schneps
\cite{Schneps:Dessins}, which contains a survey of G. Jones and D. Singmaster
on group theoretical aspects. In the graph theory literature, the original
source often cited is the 1891 article of L. Heffter
\cite{Heffter:GraphsOnSurfaces}.

\section{From a link diagram to the Kauffman state ribbon graph}
\label{A-dessin}

\begin{figure}[htbp]
   \centering
   \includegraphics[width=2.5in]{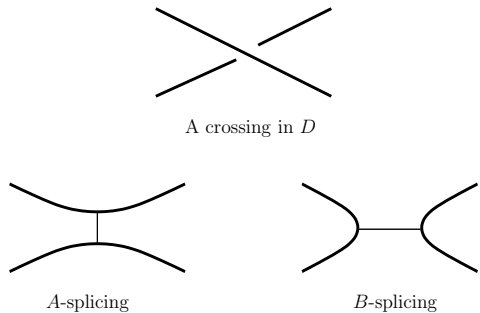} 
   \caption{Splicings of a crossing, $A$-graph and $B$-graph.}
\label{fig:AB-splicing}
\end{figure}

We will associate an oriented ribbon graph with each Kauffman state of a plane (connected) link
diagram.  The graph is constructed as follows: Given a
link diagram $D(K)$ of a knot $K$ we have, as in Figure \ref{fig:AB-splicing},
an $A$-splicing and a $B$-splicing at every crossing. For any state assignment
of an $A$ or $B$ at each crossing we obtain a collection of non-intersecting
circles in the plane, together with embedded arcs that record the crossing
splice.  Again, Figure \ref{fig:AB-splicing} shows this situation locally. In
particular, we will consider the state where all splicings are $A$-splicings. 

\medskip

To define the desired oriented ribbon graph associated to a plane link diagram, we need to define an
orientation on each of the circles resulting from the $A$ or $B$ splicings,
according to a given state assignment. (We note that these circles will become
the vertices of our ribbon graph.) For a related situation, see Vogel's
algorithm for transforming a link diagram into closed braid form
\cite{Vogel:AlexanderTheorem, BirmanBrendle:BraidSurvey}.

\begin{figure}[htbp] %
   \centering
   \includegraphics[width=3in, angle=180]{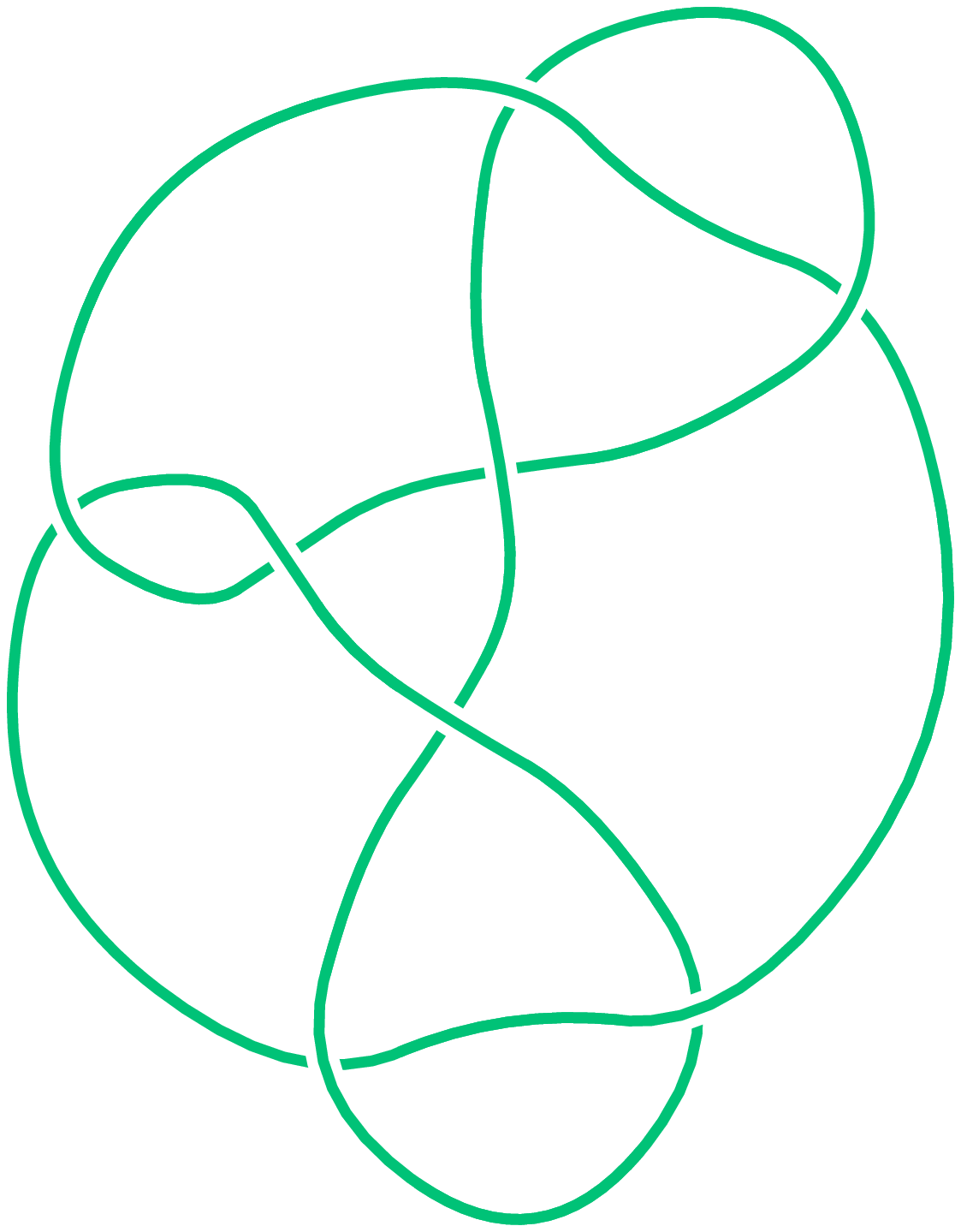} \\
   \includegraphics[width=2in]{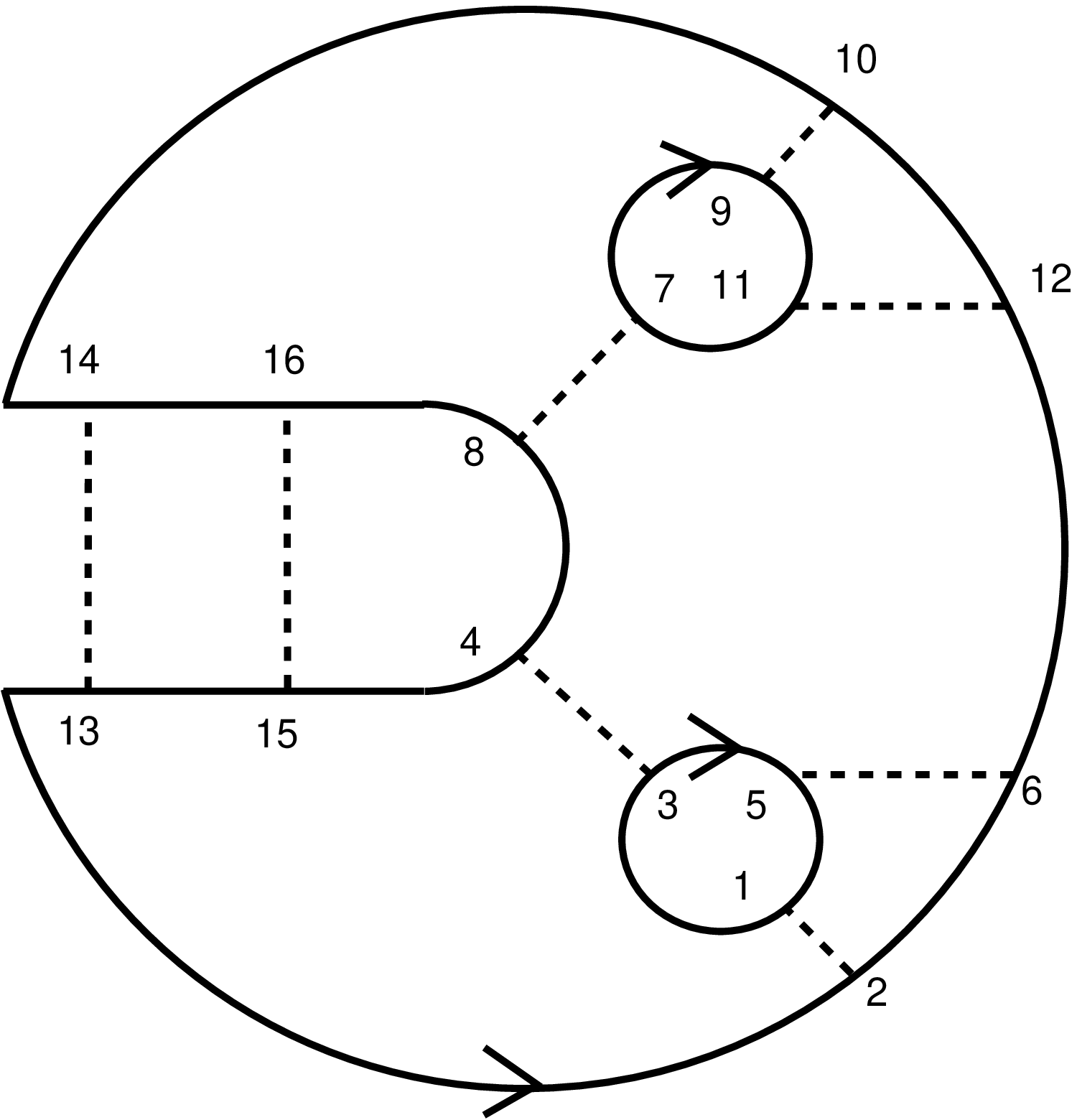} 
   \caption{The eight-crossing knot $8_{21}$ with its all-A splicing projection diagram and its orientation.}
   \label{fig:Eight21p}
\end{figure}

We orient the resulting set of circles in the plane by orienting each component clockwise or anti-clockwise according to whether the circle is inside an odd or even number of circles, respectively.
See Figure \ref{fig:Eight21p} for an example.

Given a state assignment $s:E\rightarrow\{A,B\}$ on the crossings (the eventual
edge set $E(\D)$ of the ribbon graph), the associated ribbon graph is constructed by first
resolving all the crossings according to the assigned states and then orienting
the resulting circles as above.

The set of \textit{half-edges} $\mathcal{B}$ will be the collection $\mathcal{B}=\{(C,
\gamma)\}$, where $C$ is a component circle of the resolution of the state $s$
and $\gamma$ is a directed edge from $C$ corresponding to the chosen splice at
a crossing. The permutation $\sigma_0$ permutes $\mathcal{B}$ according to the
orientation order of the endpoints of the oriented arcs $\gamma$ beginning on
$C$.  The permutation $\sigma_1$ matches the directed arc $\gamma$ with the
oppositely oriented arc beginning at the other end of the splice. We will
denote the ribbon graph associated to state $s$ by $\D(s)$.

\section{Duality}\label{sec:duality}

Given a link diagram and a state $s$, one can explicitly construct a surface
$G(s)$, following Turaev and Cromwell \cite{Cromwell:KnotsLinks, Turaev:SimpleProof}. The construction uses both 
$s$ and the \emph{dual state} $\shat$, in which every crossing is resolved in the
opposite way from $s$. It will turn out that $G(s)$ realizes the genus of $\D(s)$.

\begin{figure}[ht]
\psfrag{s}{$s$}
\psfrag{sh}{$\shat$}
\psfrag{g}{$\Gamma$}
\begin{center}
\includegraphics{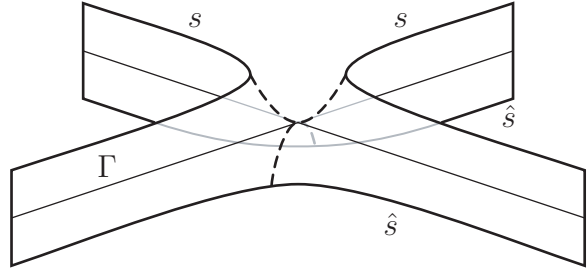}
\caption{Near each crossing of the diagram, a saddle surface interpolates 
between state circles of $s$ and state circles of $\shat$. The edges of the 
ribbon graph can be seen as gradient lines at the saddle.\label{fig:saddle}}
\end{center}
\end{figure}

Let $\Gamma \subset S^2$ be the planar, 4--valent graph of the link diagram.
Thicken the projection plane to a slab $S^2 \cross [\minus 1, 1]$, so that
$\Gamma$ lies in $S^2 \cross \{0\}$. Outside a neighborhood of the vertices
(crossings), our surface will intersect this slab in $\Gamma \cross [\minus 1,
1]$. In the neighborhood of each vertex, we insert a saddle, positioned so that
the boundary circles on $S^2 \cross \{1\}$ are the state circles of $s$, and
the boundary circles on $S^2 \cross \{\minus 1\}$ are the state circles of
$\shat$. (See Figure \ref{fig:saddle}.) Then, we cap off each state circle with
a disk, obtaining a closed surface $G(s)$. We call $G(s)$ the \emph{Turaev surface} of $s$.

\begin{lem}\label{lemma:unknotted}
$G(s)$ is an unknotted surface. In other words, $S^3 \setminus G(s)$ is a disjoint union of two handlebodies.
\end{lem}

\begin{proof}
By construction, the surface $G(s)$ has the structure of
a cell complex, whose 1--skeleton is $\Gamma$ and whose 2--cells correspond to state circles of $s$ and $\shat$. 

Thicken $\Gamma$ to a regular neighborhood $N(\Gamma)$. Because $\Gamma$ is a planar graph, $S^3 \setminus N(\Gamma)$ is a handlebody. Now, glue in the 2--cells of $G(s)$, one at a time. Each time we add a disk to the partially constructed surface, we cut the complementary manifold along a disk. It is well-known that when a handlebody is cut along a disk, the result is one or two handlebodies. Thus, by induction, $S^3 \setminus G(s)$ is a union of handlebodies.
\end{proof}

\begin{lem}\label{lemma:state-duality}
The oriented ribbon graphs $\D(s)$ and $\D(\shat)$ can both be embedded in $G(s)$. Furthermore, $\D(s)$ and $\D(\shat)$ are dual on $G(s)$: the vertices of one correspond to the faces of the other, and the edges of one correspond to the edges of the other.
\end{lem}

\begin{proof} 
Once again, we employ the crucial fact that $\Gamma$ cuts $G(s)$ into disks that correspond to state circles of $s$ and $\shat$.
These disks can be two--colored, with the $s$--disks (above $S^2 \cross \{0\}$) white and the
$\shat$--disks (below $S^2 \cross \{0\}$) shaded. 

We embed $\D(s)$ on $G(s)$ as follows. Pick a vertex in the interior of each
white $s$--disk. Then, each time two $s$--disks touch each other on opposite
sides of a crossing, connect the corresponding vertices by an edge. These edges
correspond precisely to the splicing arcs in Figure \ref{fig:AB-splicing}. 

Orienting $G(s)$ and the state
circles of $s$, in a manner compatible with the orientation of the plane projection, results in the half-edges of $\D(s)$ being  embedded with the correct cyclic ordering.

We embed $\D(\shat)$ on $G(s)$ in a similar way, by picking a vertex in the
interior of each $\shat$--disk. Now, it is easy to observe that the two ribbon graphs
are dual to each other. Every crossing of the diagram gives rise to two
intersecting edges, one in each ribbon graph. Every face of $\D(s)$ corresponds to a
shaded disk in $G(s)$, which in turn corresponds to a vertex of $\D(\shat)$ --
and vice versa. 
\end{proof}

\begin{cor}\label{cor:equal-genus}   
The genera of $G(s)$, $\D(s)$, and $\D(\shat)$ are all equal.
\end{cor}

The state surface $G(s)$ provides a concrete connection between the Jones polynomial of a link $L$ and the geometry and topology of the link complement. For example, the papers \cite{DL:VolumeIsh, FKP} use the checkerboard coloring of this surface to relate the coefficients of the Jones polynomial to the hyperbolic volume of the link. It also leads to a natural knot invariant:

\begin{df}
Given a particular projection of a link $L$, we denote the oriented ribbon graph of
the all-$A$ state by $\D(A)$, and of the all-$B$ state by $\D(B)$.
Then, we define the \emph{Turaev genus} of $L$ to be the minimum
value of $g(\D(A))$, taken over all projections of $L$. Note that by
Corollary \ref{cor:equal-genus}, this is also equal to the minimum
value of $g(\D(B))$.
\end{df}

\begin{lem}\label{lemma:alt-on-surface}
When $s$ is the all-$A$ state or the all-$B$ state, the link has an alternating projection to $G(s)$.
\end{lem}

\begin{figure}
\psfrag{sa}{$s_A$}
\psfrag{sb}{$s_B$}
\psfrag{da}{$\Downarrow$}
\begin{center}
\includegraphics[width=5in]{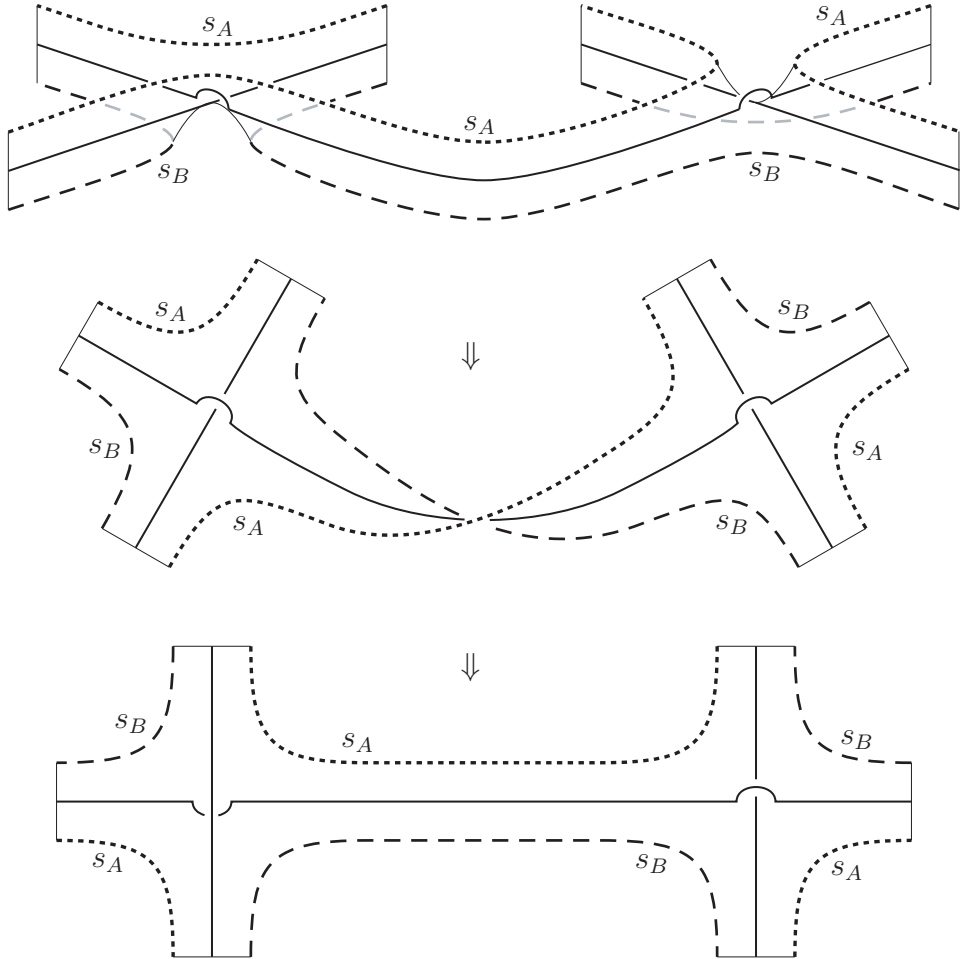}
\caption{Top: a section of the surface $G(s)$ between two consecutive over--crossings. Middle: the projection of the section. Bottom: the same section of surface, laid out flat. In this flattened picture, one can see that $L$ is alternating on the surface.\label{fig:alt-on-surface}}
\end{center}
\end{figure}

\begin{proof} When we project the link to $G(s)$, the image is the same
4--valent graph $\Gamma$. Furthermore, all the state circles of the $A$ and $B$
splicings have disjoint projections to $G(s)$, since we can draw each circle
just inside the boundary of the corresponding disk. In other words, the local picture is identical to the checkerboard coloring of an alternating diagram in the plane.

This local picture is illustrated in Figure \ref{fig:alt-on-surface}. The top panel shows an arc of $L$ between two consecutive over--crossings, as well as the corresponding section of the surface $G(s)$. When this piece of the surface is laid out flat, one can see that the left crossing becomes an under--crossing, and the projection of $L$ to the surface is alternating.
\end{proof}

\begin{cor} \label{dessingenus0}                                                                            
A knot or link $L$ has Turaev genus $0$ if and only if it is                                                
alternating.                                                                                                
\end{cor}                                                                                                   
                                                                                                            
\begin{proof} In an alternating diagram, the state circles of the $A$                                       
and $B$ splicings correspond to the checkerboard coloring of the plane. In                                          
other words, the construction of $G(s)$ recovers the (compactified) projection                                           
plane, and $\D(A)$ and $\D(B)$ have genus $0$. Conversely, if $L$ has a diagram                                        
in which  the genus of $\D(A)$ is $0$, then by Lemma \ref{lemma:alt-on-surface} the link has an alternating projection to a sphere.                                                                                                     
\end{proof}

\begin{example} \label{Example821}
 Figure \ref{fig:Eight21p} shows the non-alternating 8-crossing knot $8_{21}$, as drawn by Knotscape \cite{HTW:KnotTabulation}, and 
Figure \ref{fig:Eight21D}   the
all-$A$ associated oriented ribbon graph.

\begin{figure}[htbp] %
   \centering
   \includegraphics[width=2in]{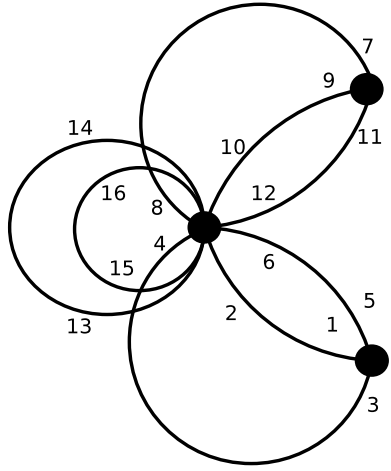} 
   \caption{All-$A$ splicing ribbon graph for $8_{21}$.}
   \label{fig:Eight21D}
\end{figure}

\begin{figure}[htbp] %
   \centering
   \includegraphics[width=3in]{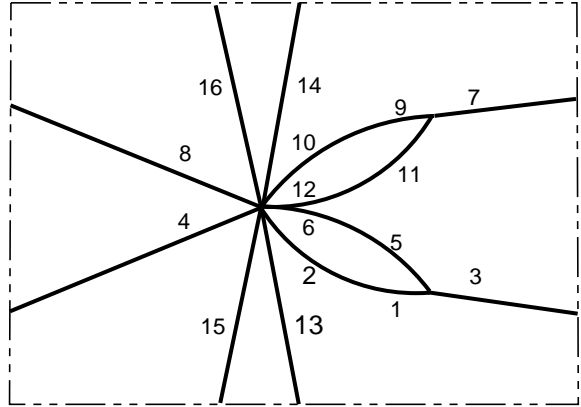} 
   \caption{Toroidal embedding of the all-$A$ splicing ribbon graph for $8_{21}$}
   \label{fig:Eight21T}
\end{figure}

With the numbering of the \textit{half-edges} as given in the diagram, the
fixed-point-free involution $\sigma_1$ is given in cycle notation by: \[
\sigma_1=\{\{1, 2\}, \{3, 4\}, \{5, 6\}, \{7, 8\}, \{9, 10\}, \{11, 12\}, \{13,
14\}, \{15, 16\}\}\] The vertex permutation reads the
half-edges around a circle of the state splicing.  In cycle notation the
permutation is:
 \[\sigma_0=\{\{2, 6, 12, 10, 14, 16, 8, 4, 15, 13\}, \{1, 3, 5\}, \{7, 9, 11\}\}.\]

By the defining property relating all three permutations, \[\sigma_2=\{\{13,
10, 7, 16, 4, 1\}, \{5, 2\}, \{8, 11, 6, 3\}, \{12, 9\}, \{15, 14\}\}\] and we
have five faces. The Euler characteristic of this ribbon graph is
\[v(\D)-e(\D) +f(\D) =3-8+5=0,\] \noindent so we have that the genus is 1, as displayed in 
\ref{fig:Eight21T}.
Since the knot is non-alternating we know by Corollary \ref{dessingenus0} that the
Turaev genus of the knot must be one.
\end{example}
\medskip

\begin{remark}
In a series of articles beginning with \cite{M98} and his book \cite{Mb2004}, Manturov introduces the concept of \textit{atom} into knot theory, motivated by Fomenko's use in his study \cite{F1991} of Hamiltonian systems.

An \textit{atom} is a triple $(M^2,\Gamma,c)$ with
\begin{enumerate}
\item $M^2$ is a closed, compact, oriented two-manifold
\item $\Gamma$ is a graph embedded in $M$ so that the complement is a disjoint union of 2-disks
\item $c$ is a two-coloring (checkerboard coloring) of the complementary regions.
\end{enumerate}

A \textit{vertical atom} is an atom where $M$ is embedded in three-space, the projection, $\pi$, to the $z$-axis is a Morse function and the level set at $z=c$, $\pi^{-1}(c)$ is a plane 4-valent graph whose complement is a disjoint union of 2-disks.  The closed surface, $G(s)$, constructed in this section is a vertical atom in which $\Gamma$ is the initial plane projection of the link (with the checkerboard coloring of the complement). And, as detailed above, the ribbon graph associated to any state, $\mathbb{D}(s)$, is embedded in this surface. In particular the graph of the all-$A$ state ribbon graph, $\mathbb{D}(A)$ can be  embedded so that the vertices are in the complementary disks of a fixed, chosen color and the edges connect the two regions of that color at a crossing.
\end{remark}

\section{The \BRT polynomial}
\label{sec:BRT}

In this section we will recall the definition of the \BRT polynomial $C(\D) \in
Z[X,Y,Z]$ of a ribbon graph $\D$ from \cite{BollobasRiordan:CyclicGraphs}. The
definition requires several different combinatorial measurements of the ribbon graph.

\begin{df}\label{def:dessin-counts}
For an oriented ribbon graph $\D$, we define the following quantities:
\begin{eqnarray*}
v(\D) &=& \mbox{the number of vertices} \; = \;  \mbox{the number of orbits of } \sigma_0,\\
e(\D) &=& \mbox{the number of edges} \; = \; \mbox{the number of orbits of } \sigma_1,\\
f(\D) &=& \mbox{the number of faces}\; = \; \mbox{the number of orbits of } \sigma_2,\\
k(\D) &=& \mbox{the number of connected components of } \D,\\
g(\D) &=& \frac{2k(\D)-v(\D) + e(\D) - f(\D)}{2}, \mbox{the \emph{genus} of } \D,\\
n(\D) &=& e(\D) - v(\D) + k(\D), \mbox{the \emph{nullity} of } \D.
\end{eqnarray*}
\end{df}

\medskip

The construction will be completed in two stages. For a ribbon graph with an edge
that is not a loop, i.e. does not connect a vertex with itself, the polynomial
satisfies some contraction/deletion relation. This reduces the computation to
the computation of the \BRT polynomial for ribbon graphs with one vertex:

\subsection{The \BRT polynomial for ribbon graphs with one vertex}\label{sec:one-vertex}

For a ribbon graph $\D$ with one vertex and $e(\D)$ edges, we define the \BRT
polynomial $C(\D)$ for one-vertex ribbon graphs as:

\[C(\D) := \sum_{\H \subset \D} Y^{n(\H)} Z^{g(\H)}.\]

Here, the summation is over all $2^{e(\D)}$ subgraphs of $\D$ obtained by
deleting a subset of edges. These are the spanning subgraphs on the same
vertex set as $\D$. %

\subsection{The \BRT polynomial of a ribbon graph with many vertices} 

Given an edge $e$ in an oriented ribbon graph, there is a naturally defined ribbon graph $\D/e$
obtained by contracting the edge $e$ to a vertex, with the cyclic order of the
half-edges now meeting at that vertex given by amalgamating the two cyclic
orders as shown below.

\begin{figure}[htbp] %
   \centering
   \includegraphics[scale=0.6]{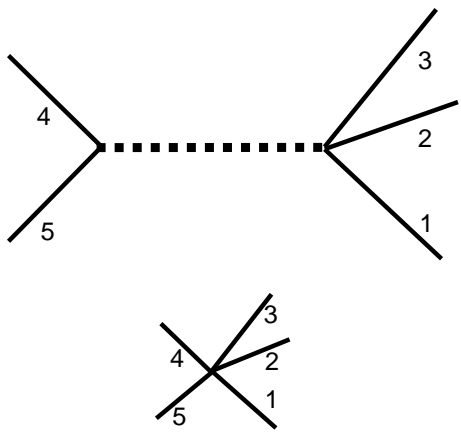} 
   \caption{Edge contraction in a ribbon graph.}
   \label{fig:Contraction}
\end{figure}

Similarly, denote by $\D-e$ the oriented ribbon graph obtained by deleting an edge $e$ and
omitting both half-edges in the orbit $e$ from the cyclic order at the
corresponding vertices. 

Recall that an edge of a graph is called a bridge when its deletion increases
the number of components by 1.

\begin{thm}[Bollob\'as-Riordan \cite{BollobasRiordan:CyclicGraphs}]  \label{BRTpolynomial} 
There is a
well-defined invariant of an oriented ribbon graph $\D$, $C(\D;X,Y,Z) \in \Z[X,Y,Z]$, the \BRT
polynomial, satisfying:

\begin{eqnarray*}
C(\D)&=& C(\D-e) + C(\D/e) \qquad \mbox{for $e$ neither a bridge nor a loop}\\
C(\D)&=& X \, C(\D/e) \qquad \qquad \qquad \mbox{for $e$ a bridge}\\
C(\D)&=&   \sum_{\H \subset \D} Y^{n(\H)} Z^{g(\H)} \qquad \qquad \mbox{if $ \D$ has one vertex}
\end{eqnarray*}

\end{thm}

\begin{remark} Note that our convention assigns the variable $X$ to a bridge
and $1+Y$ to a loop, following that of
Bollob\'as-Riordan \cite{BollobasRiordan:NonOrientableSurfaces}. The usual
convention for the Tutte polynomial assigns $X$ and $Y$, respectively.
\end{remark}

\subsection {The \BRT polynomial, the Kauffman bracket and the Jones polynomial}

Given a link projection, the Kauffman bracket  polynomial  $\langle P \rangle 
\in \Z[A,A^{-1}]$ is a regular isotopy invariant satisfying:
\begin{itemize}
\item (Normalization)  $\langle \bigcirc\rangle = 1$
\item (Trivial Components) $\langle L\coprod\bigcirc\rangle =(-A^2-A^{-2}) \langle L\rangle:=\delta 
\langle L\rangle$
\item (Skein)  $\langle L\rangle = A \langle L_A\rangle+A^{-1} \langle L_B\rangle$.
\end{itemize}

In condition three, $L_A$ is the $A$-splicing of the diagram at a chosen crossing and $L_B$ is the $B$-splicing (see Figure \ref{fig:AB-splicing}).

Recall that for alternating knots or links the Kauffman bracket of an alternating projection can be considered as an evaluation
of the Tutte polynomial of one of the checkerboard graphs, up to multiplication with a power of the variable of the Kauffman bracket  
(see \cite{Thistlethwaite:SpanningTreeExpansion}).  Chmutov and Pak have used
the \BRT polynomial to compute the Kauffman bracket of an alternating link embedded in a cylinder
$\mathbb{D}\times I$ over a ribbon graph $\mathbb{D}$.

The following theorem is a generalization of Thistlethwaite's \cite{Thistlethwaite:SpanningTreeExpansion} result to all links, not
just alternating ones. 
It does not, however, require the machinery of signed graphs needed in
\cite{Thistlethwaite:SpanningTreeExpansion} for the general case. The Kauffman
bracket is computed as an evaluation of the \BRT polynomial for a single state!

\begin{thm}\label{thm:main}
Let $\langle P \rangle \in \Z[A,A^{-1}]$ be the Kauffman bracket of a connected link projection diagram $P$ and $\D$ be the
oriented ribbon graph of $P$ associated to the all-$A$-splicing.
Then the \BRT polynomial $C(\D; X,Y,Z)$ and the Kauffman bracket are related by
\[A^{-e(\D)} \langle P \rangle = A^{2-2v(\D)} C(\D; -A^4, A^{-2} \delta, \delta^{-2}),
\]
where $\delta:=(-A^2-A^{-2})$, and $e(\D)$ and $v(\D)$ are as in Definition \ref{def:dessin-counts}.
\end{thm}

\begin{proof} 
We prove the theorem by induction on the number of crossings of $P$, or equivalently, the number of edges of $\D$. The base case is when $P$ is a single unknotted circle, and thus $\D$ has one vertex and no edges. Then it immediately follows that $A^0 \langle P \rangle =A^0 =1$. 

For the inductive step, we consider three different cases. By an abuse of notation, we will use
$\langle \D\rangle$ to denote the Kauffman bracket of a projection that yields $\D$ as its all-$A$-graph.

\smallskip
\noindent {\it \bf Case 1: $\D$ has one vertex, and thus all edges are loops.}  Then use the state sum formula for the Kauffman bracket to compute that 
{\setlength{\jot}{1.3ex}
\begin{eqnarray*}
\langle \D \rangle &=& \sum_{\H\subset\D} \delta^{f(\H)-1} A^{e(\D)  -
e(\H)}A^{-e(\H)}\\
&=& A^{e(\D)} \sum_{\H\subset\D} A^{-2e(\H)} \delta^{f(\H)-1}\\
&=& A^{e(\D)}\sum_{\H\subset\D} (A^{-2}\delta)^{2g(\H)+f(\H)-1}
(\delta^{-2})^{g(\H)}\\
&=& A^{e(\D)}\sum_{\H\subset\D} (A^{-2}\delta)^{n(\H)}
(\delta^{-2})^{g(\H)},
\end{eqnarray*}}
\noindent which is the desired specialization of the \BRT polynomial. Note that we did not need the inductive hypothesis in this case. 

\smallskip
\noindent {\it \bf Case 2: Some edge $e$ of $\D$ is a bridge.} Then
the crossing $c_e$ of $P$ corresponding to $e$ is nugatory; i.e. there is a simple
closed curve in the projection plane that intersects $P$ exactly once at the
double point $c_e$. If we smooth  the link at the nugatory crossing, the all-$A$-graph of the resulting diagram will be $\D/e$.
The properties of the Kauffman bracket imply that
$\langle \D\rangle= (-A^3) \langle \D/e\rangle$. Using this and the
inductive hypothesis,
we have
\begin{eqnarray*}
A^{-e(\D)}\langle \D\rangle &=&(-A^2)(A^{-(e(\D)-1)} \langle \D/e\rangle)\\
&=& (-A^2)(A^{2-2(v(\D)-1)}) C(\D/e;X,Y,Z)\\
&=&A^{2-2v(\D)} (-A^4)C(\D/e;X,Y,Z)\\
&=& A^{2-2v(\D)}C(\D;X,Y,Z),
\end{eqnarray*}
\noindent where the last equation comes from the second axiom of the \BRT
polynomial.

\smallskip
\noindent {\it \bf Case 3: Some edge $e$ of $\D$ is neither a bridge nor a loop.} 
Let $c_e$ be the crossing of $P$ corresponding to $e$. Note that
the all-$A$-graph of the result of an
$A$-splicing at $c_e$ will be the ribbon graph $\D-e$; and the
all-$A$-graph after performing the $B$-splicing at $c_e$ will be the ribbon graph $\D/e$. Since the projections obtained from $P$ after
each of these two splicings are connected, the induction hypothesis
applies. Using the skein relation for the Kauffman bracket, we have
\begin{eqnarray*}
A^{-e(\D)}\langle \D\rangle &=& A^{-1} (A^{-e(\D)}\langle \D/e\rangle) +
A (A^{-e(\D)}\langle \D-e\rangle) \\
 &=& A^{-2}(A^{-(e(\D)-1)} \langle \D/e\rangle) + (A^{-(e(\D)-1)}) \langle
\D-e\rangle \\
&=& A^{-2} (A^{2-2(v(\D)-1)}) C(\D/e;X,Y,Z) + A^{2-2v(\D)} C(\D-e;X,Y,Z)\\
&=& A^{2-2v(\D)}( C(\D/e;X,Y,Z) +  C(\D-e;X,Y,Z))\\
&=& A^{2-2v(\D)}C(\D;X,Y,Z),
\end{eqnarray*}
\noindent where the last equation comes from the first axiom of the \BRT
polynomial. This completes the proof of the theorem.
\end{proof}

\begin{remark} \label{rmk:AB-replacement} It is not hard to see that the
$B$-graph of a link projection $P$ is equal to the $A$-graph of the mirror image of
$P$ (compare Figures \ref{fig:AB-splicing} and \ref{fig:saddle}). It follows
that, in Theorem \ref{thm:main}, replacement of $\D(A)$ by $\D(B)$ will compute
the Kauffman bracket of the mirror image of the knot.
\end{remark}

\begin{remark}
Let $P$ be a connected projection of a link $L$, and let $w(P)$ denote the writhe of $P$.
Recall that the Jones polynomial $J_L(t)$ is obtained from
$(-A)^{-3w(P)} \langle P \rangle$ by substituting $A:=t^{-1/4}$.
Thus, by Theorem \ref{thm:main},  the Jones polynomial of $L$ is obtained as a specialization of the \BRT
polynomial of the $A$-graph corresponding to $P$.
\end{remark}

\section{The Spanning Sub-graph and Tree Expansions}
\label{sec:SpanningTrees}

The \BRT polynomial has a spanning subgraph expansion and a spanning tree
expansion, yielding the following corollaries.

For the spanning tree expansion we need an order $\prec$ of the edges of the
oriented ribbon graph $\D$. A \textit{spanning tree} $T$ of $\D$ is a subgraph with the same
vertex set, which is connected and has no cycles (equivalently zero nullity, or
no homology). For an edge $e$ of $T$ the \textit{cut} determined by $e$ and $T$
is the set of edges in $\D$ connecting one component of $T-e$ to the other. An
edge $E$ is called \textit{internally active} if $e$ is the smallest element of
the cut in the prescribed order $\prec$.  An edge $e$ not in $T$ is
\textit{externally active} if it is the smallest element in the unique cycle in
$T\cup e$.  Denote by $i(T)$ the number of internally active edges in the spanning tree $T$.

By \cite{BollobasRiordan:CyclicGraphs}  %
the spanning tree expansion of the \BRT polynomial is given by:

\[  \sum_T X^{i(T)} \sum_{S\subset\mathcal{E}(T)} Y^{n(T\cup S)} Z^{g(T\cup S)}
\]

\noindent where $\mathcal{E}(T)$ is the set of externally active edges for a
given tree $T$ (and the order $\prec$). Therefore by Theorem \ref{thm:main}, we
have:

\begin{cor}
Let $\langle P \rangle \in \Z[A,A^{-1}]$ be the Kauffman bracket of a connected link projection diagram $P$ and $\D := \D(A)$ be the oriented ribbon graph of $P$ associated to the all-$A$-splicing.
The Kauffman bracket can be computed using the fixed edge order $\prec$) by the following spanning tree expansion:
\[A^{-e(\D)} \langle P \rangle = A^{2-2v(\D)}  \sum_T X^{i(T)} \sum_{S\subset\mathcal{E}(T)} Y^{n(T\cup S)} Z^{g(T\cup S)}
\]
\noindent under the following specialization:  $\{X\rightarrow -A^4,Y\rightarrow A^{-2} \delta,Z\rightarrow  \delta^{-2}\}$, 
where $\delta:=(-A^2-A^{-2})$. \end{cor} 

The following spanning subgraph (subgraph with the same vertex set)
expansion is obtained by specializing the expansion defined in
\cite{BollobasRiordan:NonOrientableSurfaces} to the ribbon graph case where all the
edges are orientable. Again, by Theorem \ref{thm:main}:

\begin{cor}\label{prop:subgraph}
Let $\langle P \rangle \in \Z[A,A^{-1}]$ be the Kauffman bracket of a connected link projection diagram $P$ and $\D := \D(A)$ be the oriented ribbon graph of $P$ associated to the all-$A$-splicing.
The Kauffman bracket can be computed by the following spanning subgraph $\H$ expansion:
\[A^{-e(\D)} \langle P \rangle = A^{2-2v(\D)} (X-1)^{-k(\D)}\sum_{\H \subset \D} (X-1)^{k(\H)}Y^{n(\H)} Z^{g(\H)}
\]
\noindent under the following specialization:  $\{X\rightarrow -A^4,Y\rightarrow A^{-2} \delta,Z\rightarrow  \delta^{-2}\}$
where $\delta:=(-A^2-A^{-2})$.
\end{cor}

\section{Span of the Polynomial of Adequate Knots}
\label{sec:Adequate}

A connected link projection $P$ is called $A$-adequate (resp. $B$-adequate) if
and only if $\D:=\D(A)$ (resp. $\D^{*}:=\D(B)$) contains no loops (edges with
both endpoints at the same vertex). Now $P$ is called adequate if it is both
$A$- and $B$-adequate (\cite{Cromwell:KnotsLinks}), and a link is called
adequate if it admits an adequate projection. The class of adequate links
contains that of alternating ones, but it is much more general. Let $e(\D)$
denote the crossing number of $P$ and let $v(\D)$, $v'(\D)$ denote the numbers
of vertices of $\D(A)$, $\D(B)$, respectively. It is known that the span of the
Kauffman bracket of an adequate projection $P$ is given by ${\rm span} \langle
P \rangle=2e(\D) +2v(\D)+2v'(\D)-4$.  Next, we show how to derive this from the
subgraph expansion of Proposition \ref{prop:subgraph} for connected link
projections. 

The following result was essentially obtained by Kauffman \cite{Kauffman:StateModels, MortonBae:ExtremeTerms}.
Manturov \cite{M2005} put it into the atomic context described in section \ref{sec:duality}.
Here we give a proof using ribbon graphs, together with corollaries inspired by the ribbon graph approach. 
 
\begin{lem} \label{span} For a connected link projection $P$, let $M(P)$ and $m(P)$ denote
the maximum and minimum powers of $A$ that occur in $\langle P \rangle$. We have:
 
\begin{enumerate}
\item[(a)] $M(P)\leq e(\D)+2v(\D)-2$, with equality if $P$ is $A$-adequate.

\item[(b)] $m(P) \geq -e(\D)-2v'(\D)+2$, with equality if $P$ is $B$-adequate.
 \end{enumerate}
 
\noindent In particular, if $P$ is adequate then 
$${\rm span} \langle P \rangle=M(P)-m(P)=2e(\D) +2v(\D)+2v'(\D)-4.$$
\end{lem}
 
\begin{proof} Let ${\D}:={\D}(A)$ denote the oriented ribbon graph corresponding to the all-$A$
splicing of $P$.  Given a spanning subgraph $\H \subset \D$ we have, by
Definition \ref{def:dessin-counts}, $$2g(\H)= 2k(\H)-f(\H) - v(\H)+ e(\H).$$
Now a straightforward computation, using $k(\D)=1$, shows that after the
substitutions for $X,Y,Z$ given in Proposition \ref{prop:subgraph} we obtain

\begin{equation}
A^{e(\D)-2v(\D)+2}  (X-1)^{k(\H)-k(\D)}    Y^{n(\H )}  Z^{g(\H )} \; 
\rightarrow \;
A^{e(\D)-2e(\H)}   (-A^2 -A^{-2})^{f(\H)-1}. 
\end{equation}

Then $M(\H):=e(\D)-2 e( \H)+2 f(\H) -2$ is the highest power of $A$ contributed
by $\H$. Let $\H_0$ denote the spanning subgraph that contains only the
vertices of $\D$ and no edges. Now every spanning subgraph $\H \subset \D$
is obtained from $\H_0$ by adding a number of edges. This can be done in stages
so that there are subgraphs $\H_0, \ldots , \H_k$, with $\H_k=\H$ and such
that, for $i=1,\ldots, k$, $\H_i$ is obtained from $\H_{i-1}$ by adding exactly
one edge. Then, $e(\H_{i})=e(\H_{i-1})+1$ and $f(\H_{i})=f(\H_{i-1})\pm 1$, and
hence the difference $M(\H_{i-1})-M(\H_{i})$ is 0 or 4.  Thus we obtain
\begin{equation} \label{equation72}
M(\H)\leq M(\H_{0})= e(\D)+2v(\D)-2,
\end{equation} 
for every spanning subgraph $\H \subset \D$. Now if $P$ is $A$-adequate then,
since an edge is added between two distinct vertices of $\H_0$, we have
$f(\H_{1})=f(\H_{0})- 1$, and so $ M(\H_{i})$ decreases at the first step. By
(\ref{equation72}) it never increases, so we have

\begin{equation}
M(\H)< M(\H_{0})= e(\D)+2v(\D)-2, 
\end{equation}

\noindent
for every $\H\neq \H_0$. Thus the term with degree $M(\H_{0})$ is never
canceled in the subgraph expansion of $\langle P\rangle$ and part (a) of the
lemma follows. Part (b) follows by applying part (a) to the mirror image of $P$
and using the observation in Remark \ref{rmk:AB-replacement}.
\end{proof}
 
\begin{cor} \label{corspan}    Let $P$ be an adequate projection of a link $L$ and let $\D(A)$ and
$\D(B)$ be as above. Then the genus $g(\D(A))=g(\D(B))$ is an invariant of the
link $L$.
\end{cor}
 
\begin{proof} 
It is known that $e(\D)$
is actually the minimal crossing number of the link $L$ and thus an invariant
of $L$ \cite[\S 9.5]{Cromwell:KnotsLinks}. It is also known that the span of
the Kauffman bracket of any link projection is an invariant of the link. Since
$P$ is adequate, by Lemma \ref{span}, we have ${\rm span} \langle
P\rangle=2e(\D) +2v(\D)+2v'(\D)-4$. By Lemma \ref{lemma:state-duality}, we have
$v'=f(\D(A))$ and thus (Definition \ref{def:dessin-counts}) $2-2
g(\D(A))=v(\D)+v'(\D)-e(\D)$. Thus $4g(\D(A))=4e(\D) -{\rm span} \langle
P\rangle$, and the conclusion follows from the previous observations.
\end{proof}

The span of the Kauffman bracket  in the variable $A$  of a knot projection is four times the span of the Jones polynomial in the variable $t$.

\begin{cor} \label{genusestimate} Let $P$ be a $c$-crossing, connected projection of a link $L$,
let $s_L$ denote the span of the Jones polynomial of $L$ and let $g_L$ be the Turaev genus of the link.
Then,
$$g_L \leq c- s_L.$$
\end{cor}

\begin{proof}
By Lemma \ref{span} and its proof, we have that the span of the Kauffman bracket is less or equal to $4c-4g(\D(A))$. Since the span of the Kauffman bracket  in the variable $A$  of a knot projection is four times the span of the Jones polynomial in the variable $t$, we have  $g_L\leq g(\D(A)) \leq c- s_L$.
\end{proof}

\begin{remark}
The estimate in Corollary \ref{genusestimate} is sharp for some families but not in general. For example, it is sharp for the family of non-alternating pretzel knots  $P(a_1, \ldots, a_r, b_1, \ldots, b_s)$, where $a_i\geq 2$, $b_j\geq 2$, $r,s \geq 2$, as in Figure \ref{fig:pretzel}. On the one hand, these knots are non-alternating, so they have Turaev genus at least one. On the other hand, Lickorish and Thistlethwaite \cite{LT} prove that the span of the Jones polynomial is one less than the crossing number. Thus the estimate of Corollary \ref{genusestimate} is sharp.

\begin{figure}[htbp]
   \centering
   \includegraphics[width=3.5in]{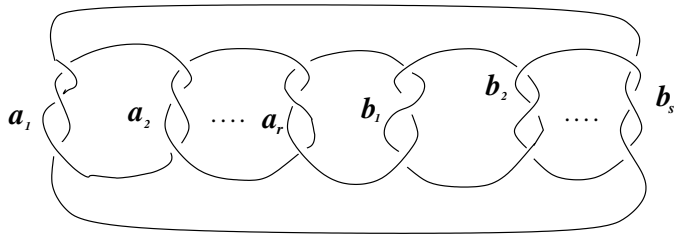} 
   \caption{Non-alternating pretzel knots that realize the estimate of Corollary \ref{genusestimate}.}
\label{fig:pretzel}
\end{figure}

For a non-sharp example, consider the $8$-crossing knot $8_{21}$ in Figure \ref{fig:Eight21p}.
The span of its Jones polynomial polynomial is $6$ (Knotscape  \cite{HTW:KnotTabulation}), and its Turaev genus is one. 
\end{remark}

\bibliography{BRTv4.bbl}
\bibliographystyle {alpha}
 
\end{document}